\documentclass[reqno, 11pt]{amsart}
\usepackage{amssymb,latexsym, mathabx}
\usepackage{enumerate,subcaption,color,bbm}
\usepackage[utf8]{inputenc}

\usepackage{url}
\usepackage{mathrsfs}
\makeatletter
\@namedef{subjclassname@2010}{%
  \textup{2010} Mathematics Subject Classification}
\makeatother

\renewcommand{\eqref}[1]{(\ref{#1})}

\newtheorem{thm}{Theorem}[section]

\newtheorem{lem}[thm]{Lemma}

\newtheorem*{cor}{Corollary}

\def\1{\mathbbm{1}}

\theoremstyle{definition}
\newtheorem{remark}[thm]{Remark}

\usepackage{pgf,tikz,pgfplots}

\usepackage{mathrsfs}
\newcommand{\image}[2]{\begin{center}
\includegraphics[width=#1\textwidth]{#2.png}
\end{center}}

\numberwithin{equation}{section}

\DeclareMathOperator{\Ocal}{\mathcal{O}}
\DeclareMathOperator{\sym}{sym}

\newcommand{\myF}{\mathscr{F}}
\newcommand{\myG}{\mathscr{G}}
\newcommand{\myH}{\mathscr{H}}

\newcommand{\R}{{\mathbb R}}
\newcommand{\Q}{{\mathbb Q}}
\newcommand{\N}{{\mathbb N}}
\newcommand{\abs}[1]{\left|#1\right|}
\newcommand{\ve}{\varepsilon}

\begin{document}

\title[A Modular Analogue]{A Modular Analogue of a Problem of  Vinogradov}

\author{R. Acharya}
\email{ratnadeepacharya87@gmail.com.}
\address{Ramakrishna Mission Vivekananda Educational Research Institute,
P.O. Belur Math, Dist. Howrah
West Bengal, India.
PIN: 711 202.} 

\author{S. Drappeau}
\email{sary-aurelien.drappeau@univ-amu.fr}
\address{Institut de Math\'ematiques de Marseille, U.M.R. 7373 \\
Aix Marseille Universit\'e, Centrale Marseille, CNRS \\
Campus de Luminy, Case 907 \\
13288 Marseille Cedex 9, France}
\author{S. Ganguly}
\email{sgisical@gmail.com}
\address{Theoretical Statistics and Mathematics Unit \\
Indian Statistical Institute \\
203 Barrackpore Trunk Road, Kolkata 700108, India}
\author{O. Ramar\'e}
\email{olivier.ramare@univ-amu.fr}
\address{CNRS / Institut de Math\'ematiques de Marseille \\
Aix Marseille Universit\'e, U.M.R. 7373 \\
Site Sud, Campus de Luminy, Case 907 \\
13288 MARSEILLE Cedex 9, France}

\today
\begin{abstract}
  Given a primitive, non-CM, holomorphic cusp form~$f$ with normalized Fourier
  coefficients~$a(n)$ and  given an interval~$I\subset [-2, 2]$, we study the least prime~$p$ such
  that~$a(p)\in I$ . This can be
  viewed as a modular form analogue of Vinogradov's problem on the least
  quadratic non-residue. We obtain strong explicit bounds on~$p$,
  depending on the analytic conductor of~$f$ for some specific
  choices of~$I$.
\end{abstract}

\maketitle

\section{Introduction}

The present article is concerned with understanding the 
distribution of the initial Fourier coefficients of primitive holomorphic cusp forms
at primes. Suppose $f$ is such a  form  of weight $k$
for the group $\Gamma_0(N)$. We further assume that $f$ is non-CM and has trivial nebentypus. 
The normalized Fourier coefficients of~$f$
at infinity are denoted by~$(a(n))_{n\geq 1}$,  so
that~$a(1) = 1$ and 
\[
 f(z)=\sum_{n=1}^{\infty} a(n)n^{\frac {k-1}{2}}e(nz),
\]
where, as usual, $e(z)$ denotes $e^{2 \pi iz}$ and with this
normalization, the Ramanujan bound (proved by Deligne \cite{Deligne1974})
says $-2\leq a(p)\leq 2$ for primes~$p$. Furthermore, the function~$n\mapsto a(n)$ is real-valued and
multiplicative. We refer the reader to the text~\cite{Iwaniec1997} for background information on holomorphic modular forms.
The Sato-Tate conjecture for distribution of the angles $\theta_p$, defined by $a(p)=2\cos \theta_p$, as $p$ runs over primes,
which is now a theorem of Clozel, Harris, Shepherd-Barron and
Taylor~\cite{Clozel-Harris-Taylor*08,Taylor*08,HarrisEtAl2010}, implies, in
particular, that any interval of positive measure within $[-2,2]$
contains infinitely many values of $a(p)$. The goal of this article
 is to obtain bounds for the least prime $p$ such that $a(p)$
lies in a fixed interval $I \subset[-2,2]$. This can be considered as an analogue of Vinogradov's
problem of estimating, given a modulus~$q\ge 1$, the size of the least
quadratic non-residue modulo~$q$ (see~\cite{Burgess*57},
  \cite{Vinogradov*27b}). The quality of our bounds will be measured in terms of the analytic conductor~$q(f) = Nk^2$ of the form $f$ (see \S 2.1), and also separately in term of 
  the weight $k$ of the form, considering the level $N$ to be fixed and 
  in terms of the level $N$, considering the weight $k$ to be fixed.
  We restrict our attention to forms with trivial nebentypus in order to clarify
  the presentation but the methods presented here can be extended to a more general setting.
 
Let~$I\subset[-2, 2]$.  Theorem~1.6 of the
paper ~\cite{LemkeOliver-Thorner*07} of Lemke-Oliver and Thorner
implies that there exists a constant $A$ depending only on~$I$
such that $a(p)\in I$ for some prime~$p\leq q^A$. 
Their method
relies on effective log-free zero density estimates for the
$L$-function associated with~$f$, and the Tur\'an power-sum
method. The value of the constant $A$ is not stated explicitly in their paper but it is not hard to see that the constant is effective and can be
worked out explicitly. However the method is likely to produce quite large
values of~$A$.
Our aim in the present work is to make the value of  $A$ as small as possible for some specific intervals.

We define, when $\kappa$ is positive and $x\in[0,1]$: 
\begin{equation}
  \label{defmyF}
  \myF(x;\kappa)=\int_0^{x/(1+x)}\frac {h^{\kappa-1}dh}{1-h}
  =\sum_{k\ge0}\frac{1}{\kappa+k}\Bigl(\frac{x}{1+x}\Bigr)^{\kappa+k}.
\end{equation}
Note that $\myF(\cdot;\kappa)$ is increasing between $\myF(0;\kappa)=0$ and $\myF(1;\kappa)=\int_0^{1/2} \frac {h^{\kappa-1}dh}{1-h}$. We thus define a function $\myG(\cdot;\kappa)$ with value in $[0,1]$ by
\begin{equation}
  \label{eq:8}
  \myG(y;\kappa)=\max\{x\in[0,1]: \myF(x;\kappa)\leq 1/y\}.
\end{equation}
The function $\myG$ is non-increasing and we have $\myG(y;\kappa)=1$ when $y\leq 1/\myF(1;\kappa)$ and by convention~$\myG(\infty;\kappa)=0$. 

We now state our main results which depend crucially on knowledge
about the analytic properties of the symmetric power $L$
functions associated to~$f$ (see \S \ref{Sympower} for definition). This is likely to change in the future;
only small changes would be required in our proofs to reflect any such
improvement. Here is the assumption we rely on.

\textbf{Hypothesis} $\mathbf{\myH_\ell}$: The $L$-function  $L(s, \sym^{\ell}(f))$ has analytic continuation to the entire complex plane  and it satisfies the bound
\[
 L(1/2+it, \sym^{\ell} (f)) \ll_{\ve} q(\sym^{\ell} (f),s)^{\lambda_{{\ell}}+\ve}
\]
 for any
$\ve>0$.

For holomorphic forms, the  automorphy of $L(s,\sym^\ell f)$  has been known for $\ell\leq 8$ by~\cite{GelbartJacquet1978,KimShahidi2002b,Kim2003,KimShahidi2002a,ClozelThorne2015,ClozelThorne2017}, and has recently been proved for \emph{all} $\ell$ when~$N$ is squarefree  by Newton-Thorne~\cite{NewtonThorne2020}.
 As a result, these $L$-functions admit holomorphic continuation to the entire complex plane and by the convexity principle, $\myH_\ell$ holds with $\lambda_\ell=1/4$ (known as the convexity bound) for $\ell \leq 8$ unconditionally and for all $\ell$ when~$N$ is squarefree.

Our results are the following. 

\begin{thm}
  \label{thm1}
  For any $\delta\in(0,2]$, let $\theta_1(\delta)=\myG(2+\delta; \delta)$. The function $\theta_1$ is increasing and we have $\theta_1(0+)=0$ and $\theta_1(1)=0.3956\cdots$. Suppose $\lambda_1>0$ is an exponent that satisfies the hypothesis $\myH_\ell$ below for $\ell=1$, and let~$\ve>0$. Then for
  $q=N$ or $k^2$ sufficiently large,
   there exists a prime 
   \[
   p\ll_\ve q^{\frac{2\lambda_1}{1+\theta_1(\delta)}+\ve}
   \]
    with $a(p)\leq \delta$.
\end{thm}

\begin{remark}
The convexity bound (Phragm\'en-Lindel\"of principle) allows taking $\lambda_1=\frac{1}{4}$ but better exponents, called subconvex exponents are known in both the weight
and the level aspects. For example, one may take $\lambda_1=\frac{1}{6}$ when $N=1$ by a result of Jutila and Motohashi \cite{JM}.
\end{remark}
\begin{thm}
  \label{thm2}
  For any $\delta\in(0,1]$, let $\theta_2(\delta)=\myG((1+\delta)^2;2\delta+\delta^2)$. The function $\theta_2$ is increasing when $\delta\leq 0.5305\cdots$, and constant equal to $1$ afterwards.  We have $\theta_2(0+)=0$, $\theta_2(1/2)=0.9093\cdots$, $\theta_2(1)=1$. Suppose $\lambda_2>0$ is an exponent that satisfies the hypothesis $\myH_\ell$ below for $\ell =2$, and let~$\ve>0$. For any  $\delta\in[0,1]$, and for $q=N$ or $k^2$ sufficiently large, there exists a prime 
  \[  
    p\ll_\ve q^{\frac{4\lambda_2}{1+\theta_2(\delta)}+\ve}
  \] 
  with $|a(p)|\leq 1+\delta$.
\end{thm}
\begin{remark}
The convexity bound allows the choice $\lambda_2=\frac{1}{4}$ and currently this is the best known exponent. Obtaining a subconvex estimate for
the symmetric square $L$-function in the level or the weight aspect is a challenging problem.
\end{remark}

It turns out that showing the existence of  primes $p$ of  small size in terms of  the conductor (i.e., weight and level) such that $a(p)\geq 0$ is rather difficult. By utilizing the fact that 
hypotheses~$\myH_\ell$ holds true for~$1\leq\ell\leq 5$, we are able to show the following result:
\begin{thm}\label{thmpos}
There is a prime~$p \ll k^{24} N^{21}$ such that~$a(p)\geq 0$.
\end{thm}

The results above are all obtained using a similar strategy and this is summarized in Theorem~\ref{geneV} below. For some specific intervals, however, we  obtain better bounds by employing ad hoc techniques using $L$-functions as we now describe.

\begin{thm} \label{thmnew}
For any $\ve>0$, there is a prime $p=\Ocal_{\ve}(kN)^{1+\ve}$ such that $a(p) <0$. 
\end{thm}

\begin{cor}
  The least prime such that $a(p)\neq 0$ is $\ll_{\ve} (kN)^{1+\ve}$, for any $\ve>0$.
\end{cor}

\begin{remark}
As the proof of the above theorem shows, the exponent $1$ can be replaced by $4\lambda_2$ and any subconvex estimate $\lambda_2 <1/4$ for the symmetric square $L$-function will lead to an improvement of the above result. 
\end{remark}

The next result relates the possibility of the initial coefficients at primes assuming extreme values with the size of $L(1,f)$. For 
$q=Nk^2$, let
$$ \gamma^- := \liminf_{q\to\infty} \frac{\log L(1, f)}{\log\log q}, \qquad \gamma^+ := \limsup_{q\to\infty} \frac{\log L(1, f)}{\log\log q}. $$
From the zero-free region of $L(s, f)$ (See \cite{HR}), the standard techniques yield
\begin{equation}
  -2\leq \gamma^- \leq \gamma^+ \leq 2.\label{eq:bounds-L1f-uncond}
\end{equation}

\begin{thm} \label{thm3}
  For any~$\delta, \ve>0$, the least prime $p$ such that $a(p)>\gamma^--\delta$ is
  $\Ocal(q^\ve)$.
  Similarly, the least prime~$p$ such that~$a(p)<\gamma^++\delta$ is~$\Ocal(q^\ve)$.
\end{thm}

\begin{remark}
The bounds~\eqref{eq:bounds-L1f-uncond} seem to be the best known, and any improvement would yield a non-trivial result in Theorem~\ref{thm3}. The quality of the upper-bound on~$p$, namely~$\Ocal(q^\ve)$, compared to the above results, suggests that improving the bounds~\eqref{eq:bounds-L1f-uncond} is a difficult task. Under the Riemann Hypothesis for $L(s, f)$, one has the bounds
\[
(\log\log q)^{-2} \ll L(1, f) \ll (\log\log q)^{2},
\]
at least in the case $N=1$ (see \cite[Thm. 3]{LW} for a precise and stronger statement), which yields conjecturally~$\gamma^- = \gamma^+ = 0$. Furthermore, it is known that these bounds hold for almost all forms (see \cite[Cor. 2]{LW2} for a precise statement). 
\end{remark}

Several authors investigated the smallest \emph{integer} $n$ such that
$a(n)<0$, see for instance \cite{Iwaniec-Kohnen-Sengupta*07},
\cite{Kowalski-Lau-Soundararajan-Wu*10}, \cite{Lau-Liu-Wu*12} or
\cite{Matomaki*12}. It follows from~\cite{Matomaki*12} that the least
such~$n$ is~$\Ocal(q^{3/8})$, where $q=Nk^2$. A closer scrutiny of their proofs reveals
that the integer~$n$ they produce is either a prime or the square of a
prime. Indeed, all the above works make use of the contrast between the sizes of $a(p)$ and $a(p^2)$ 
forced by the Hecke relation $a(p)^2-1=a(p^2)$ for primes $p$. Since we aim at localizing
only $a(p)$'s, the coefficients at primes,  we cannot rely on such procedures. In fact, the two
methods we propose are \emph{reverse}: from a localization on $a(p)$, we
show that some polynomial in $a(p)$ has to be large for many
primes~$p$. This polynomial defines the value at $p$ of a new function
whose Dirichlet series we approximate with products of $L(s,\sym^\ell f)$ and it is by using the analytic properties of these latter that
we reach a contradiction. To find an integer $n$ such that $a(n)<0$,
only the analytic properties of $L(s,f)$ are required.

Regarding  bounds conditional on the Riemann Hypothesis,
Ankeny~\cite{Ankeny*52} has proved that for any non-trivial
character~$\chi\mod{q}$, if the Riemann hypothesis is true
for~$L(s,\chi)$, then the least~$n$ such that~$\chi(n)\neq 1$
is~$\Ocal((\log q)^2)$. It is not difficult to show that the analogous
phenomenon holds in our setting:
\begin{thm}\label{th:ankeny}
  Assume that for all~$\ell\geq 1$, the function~$L(s,\sym^\ell f)$ is
  entire and satisfies the Riemann hypothesis. Then for any
  interval~$I\subseteq [-2, 2]$ of positive measure, the least
  prime~$p$ such that~$a(p)\in I$ satisfies~$p\ll_I (\log q)^2$.
\end{thm}

\bigskip

Let us now state our general theorem depending on the hypothesis $\mathbf{\myH_\ell}$. Note that this result implies Theorems~\ref{thm1}, \ref{thm2} and \ref{thmpos}.

\begin{thm}[Generic theorem]
  \label{geneV}
  Let $(b_\ell)_{1\leq \ell\leq L}$ be non-negative integers, Let $\kappa>0$ and~$F$ be real, and let $I\subset[-2,2]$ be such that
  \begin{equation}
    \label{eqI}
    \left\{ \begin{aligned} & \forall x\in [-2, 2]\smallsetminus I, \quad \sum_{1\leq \ell\leq L}b_\ell U_\ell(x/2)\geq \kappa>0, \\ & \forall x\in[-2, 2], \quad \sum_{1\leq \ell\leq L}b_\ell U_\ell(x/2)\geq F, \end{aligned}\right.
  \end{equation}
  where $U_\ell$ are the Chebyshev polynomials of the second kind. Then, on assuming $(\myH_\ell)_{\ell\leq L}$, the least prime $p$ such that $a(p)\in I$ satisfies 
  
    \begin{equation}\label{vN}
     \frac{\log p}{\log N} \leq \frac{2\sum_{\ell}\ell b_{\ell} \lambda_{\ell}}{1+\myG(\kappa-F;\kappa)}+\ve,
  \end{equation}
  for any $\ve>0$ and~$N$  large enough with respect to the weight $k$ and  $\ve$;
  and 
  \begin{equation}\label{vk}
      \frac{\log p}{\log k} \leq \frac{2\sum_{\ell}(\ell+\epsilon(\ell)) b_{\ell} \lambda_{\ell}}{1+\myG(\kappa-F;\kappa)}+\ve.
  \end{equation}
  for any $\ve>0$ and $k$  large enough with respect to the level $N$ and $\ve$. Here~$\epsilon(\ell) = \frac{1-(-1)^\ell}2 \in\{0, 1\}$ is the parity of~$\ell$.
\end{thm}

\begin{figure}
  \centering
  \begin{subfigure}{.5\textwidth}
    \image{0.9}{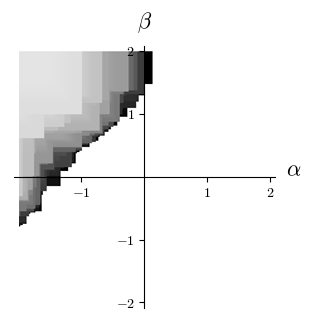}
  \end{subfigure}%
  \begin{subfigure}{.5\textwidth}
    \image{1.3}{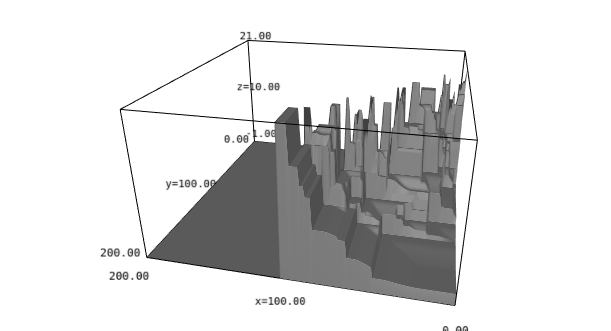}
  \end{subfigure}
  \caption{Upper-bound on~$\frac{\log p}{\log N}$ in Theorem~\ref{geneV} for~$I=[\alpha, \beta]$.}
  \label{fig:geneV}
\end{figure}

The intervals~$[\alpha, \beta]$ for which there is a linear combination
with non-negative coefficients of~$U_1, \dotsc, U_8$ which takes
positive values outside~$[\alpha, \beta]$
delimit a curve in~$(\alpha, \beta)$, whose
exact determination is an interesting question (without the
non-negativity condition, the analogue for~$U_1, \dotsc, U_4$ was
solved in Appendix~A of~\cite{LemkeOliver-Thorner*07}). Between this
curve and the diagonal~$\alpha=\beta$, Theorem~\ref{geneV} yields an
upper-bound on~$\frac{\log p}{\log q}$, which gets smaller as one
moves away from the diagonal. This is represented in
Figure~\ref{fig:geneV}, which was obtained by case-by-case analysis of
all linear combinations with~$\sum_{\ell\leq 8} \ell b_\ell \leq
42$. On the left, darker colors indicate a larger upper-bound.

Theorem~\ref{geneV} should be compared with Theorem 1.8
of~\cite{LemkeOliver-Thorner*07}. In both cases, we are given an
interval~$I\subset[-2, 2]$, and we are looking for the least prime~$p$
such that~$a(p)\in I$. In Theorem 1.8
of~\cite{LemkeOliver-Thorner*07}, the authors obtain an exponent
depending on the quality with which the indicator
function~$\1_{I}$ can be minorized by a
linear combination of~$U_0, U_1, U_2, \dotsc$. In Theorem~\ref{geneV}, we obtain an exponent depending on the quality with which the complementary indicator function~$\1_{[-2, 2]\smallsetminus I}$ is minorized by a linear combination \emph{with non-negative coefficients} of~$U_1, U_2, \dotsc$. An inconvenient of our method is that there is no clear description of the allowable intervals~$I$. 
Theorems \ref{thm1}-\ref{thmpos} indicate that, when it can be applied, the
method described here yields non-trivial numerical results.

\subsection*{Notation} Our notation is quite standard. We follow the usual practice of denoting by $p$ an arbitrary prime and by $\ve$ an arbitrarily small positive real number which need not be the same in every occurrence. For any set $X \subset \mathbb{R}$ and maps $F:X\mapsto \mathbb{C}$ and $G: X \mapsto [0, \infty)$, we write
\[F(x) \ll G(x) \ \text{or} \ F(x)={\Ocal}(G(x))\] if there exists a   $C> 0$  such that 
$|F(x)| \leqslant C G(x)\  \textnormal{ for all }\  x \in X.$ Sometimes, the implied constant $C$ depends on some parameters and this dependence is shown in the subscript. For example, often the implied constant depends on the parameter $\ve$, an arbitrarily small positive real number and we display this dependence by writing $\ll_\ve$ or $\Ocal_{\ve}$.  Sometimes, the dependence is not shown when it is clear from the context in order to avoid making the notation too cumbersome. By 
 $\check{\eta}$, we denote the
Mellin transform of a  function $\eta$:
\begin{equation}
  \label{eq:10}
  \check{\eta}(s)=\int_{0}^\infty \eta(t)t^{s-1}dt. 
\end{equation}

\section*{Acknowledgements}
The first author would like to thank Ramakrishna Mission Vivekananda Educational and Research Institute for a wonderful academic atmosphere and the National Board for Higher Mathematics (DAE file reference  no.~0204/12/2018/R \& D II/ 6463) for financial support.\\
 The third author  thanks  CNRS (Centre National de la Recherche Scientifique) 
through its Laboratoire International Associé “Indo-French Program for
Mathematics” and ISI (Indian
Statistical Institute) for financial support that enabled him to visit Aix-Marseille Université  in June, 2017, where this work was initiated. He also thanks the department of Mathematics of Aix-Marseille Université for its warm hospitality. 

\section{Background on modular forms and $L$-functions}

\subsection{Symmetric power $L$-functions}\label{Sympower}
For a primitive form $f$, as in the introduction, its normalized coefficients $a_f(p)=a(p)$ can be written as
\[
a(p)=\alpha_f(p)+\beta_f(p)
\]
where, for $p\not\mid N$, $\alpha_f(p)=1/\beta_f(p)$ and both are complex numbers of absolute value~$1$.
For each~$\ell\in\N$, the~$\ell$-th symmetric power $L$-function of~$f$ is defined, for $\Re s>1$, by
\begin{equation}
  \label{defsymk}
  L(s,\sym^\ell f)= \prod_{p} \prod_{0\leq j\leq \ell}\bigl(1-\alpha_f(p)^{\ell-j}\beta_f(p)^j/p^s\bigr)^{-1}
  =: \sum_{n\geq1}\frac{a_{\sym^\ell f}(n)}{n^s}.
\end{equation}
We have $\sym^1f=f$ and it is convenient to set $\sym^0f=\1$ so that $ L(s,\sym^0 f)=\zeta(s)$.
It is expected from a general conjecture of Langlands \cite{Langlands} that for every $\ell$, there is a cuspidal automorphic representation 
of~$GL_{\ell+1}(\mathbb{A}_\Q)$ that corresponds to the $L$-function~$L(s, \sym^\ell f)$.
For~$1\leq \ell \leq 8$, this  was shown in~\cite{GelbartJacquet1978} (for~$\ell=2$), \cite{KimShahidi2002b} (for~$\ell=3$), \cite{Kim2003,KimShahidi2002a} (for~$\ell=4$) and~\cite{ClozelThorne2015,ClozelThorne2017} (for~$5\leq \ell \leq 8$). When~$N$ is squarefree, this has been announced for all~$\ell\geq0$ in \cite{NewtonThorne2020}.

Following  \cite[Eq.(5.5)]{Iwaniec-Kowalski*04}), we define the analytic conductor of  $L(s,\sym^\ell f)$ as
\begin{equation}
  q(s,\sym^{\ell} (f))=N^{\ell}(|t|+2)^{\ell+1}k^{\ell+\epsilon(\ell)},\label{eq:def-qsymell}
\end{equation}
with~$\epsilon(\ell) = \frac{1-(-1)^\ell}2$ being $1$ or $0$ according as $\ell$ is odd or even, as in the statement of Theorem~\ref{geneV}.

 Once we know that a symmetric power $L$-function comes form an automorphic representation, the analytic continuation and functional equation for that $L$-functions follows from~\cite{GodementJacquet} and thus the 
Phragm\'en-Lindel\"of convexity principle (or the approximate functional equation~\cite[eq. (5.20)]{Iwaniec-Kowalski*04}) implies that
for~$1\leq \ell\leq 8$, the hypothesis~ $\mathbf{\myH_\ell}$ holds with the value  $\lambda_{\ell} =1/4$, even for $\delta=0$. This is known as the convexity bound. Giving
a bound on an $L$-function that is stronger than the convexity bound is a challenging problem which has been solved in a few cases (see \cite{Munshi} and the references therein) and
this is known as the subconvexity problem.
Sometimes we are interested in the size of the $L$-functions in terms of only the size of the variable $t$, or the weight $k$ or the level $N$. A result 
of Jutila and Motohashi \cite{JM} says that taking $\lambda_1 =1/6$ is permissible in the weight and the $t$-aspect. We further define 
\begin{equation}
 q(\sym^{\ell} (f)):= N^{\ell} k^{{\ell}+\epsilon({\ell})}.
\end{equation}
In particular, $q(f) =Nk^2$ and $q(\sym^2 (f))=N^2k^2$. Note that in
the weight aspect, $q(f)$ and $q(\sym^2 (f))$ are of the same order.

For the coefficients of the symmetric $\ell$-th power $L$-function of $f$, we have the
following relation for every prime $p$:
\begin{equation}
  \label{eq:5}
  a_{\sym^\ell f}(p)=a\bigl(p^\ell\bigr)=U_\ell(\cos\theta(p))=U_{\ell}(a(p)/2)
  =\frac{\sin((\ell+1)\theta(p))}{\sin\theta(p)},
\end{equation}
where $U_\ell$ is the Chebyshev polynomial of second kind, whose properties we recall next.

\subsection{Chebyshev polynomials of the second kind}

\label{Che}

We recall that the Chebyshev polynomial of second kind $(U_\ell)_{\ell\geq 0}$ are defined by
\begin{equation}
  \label{eq:6}
  U_0=1, \qquad U_1=2x, \qquad U_{\ell+1}- 2x U_\ell+U_{\ell-1}=0.
\end{equation}
These polynomials form an orthonormal basis in the space of polynomials on the interval $[-1, 1]$ relative to the Hermitian product
\begin{equation}
  \label{eq:7}
  \langle f, g \rangle = \int_{-1}^1 f(x)\overline{g(x)}\tfrac{2}{\pi}\sqrt{1-x^2}dx.
\end{equation}
The first few are given by
\begin{align*}
  U_2&=4x^2-1,\\
  U_3&=8x^3-4x,\\
  U_4&=16x^4-12x^2+1,\\
  U_5&=32x^5-32x^3+6x,\\
  U_6&=64x^6-80x^4+24x^2-1,\\
  U_7&=128x^7-192x^5+80x^3-8x,\\
  U_8&=256x^8-448x^6+240x^4-40x^2+1.
\end{align*}
The last equality in Eq. \eqref{eq:5} comes from the relation
\[
U_{n}
 (\cos \theta) =\frac{\sin((n+1)\theta)}{\sin \theta}.
 \]

\section{Auxiliary Lemmas}

\subsection{Convolutions}

\begin{lem}
  \label{genDir}
 Assume $(\myH_\ell)_{1\leq \ell\leq L}$.
   Let $L\geq 1$ be an integer and let $(b_\ell)_{0\leq \ell\leq L}$ be a collection of non-negative integers. Then, we have the equality
    \begin{equation*}
      \prod_{p}\biggl(1+\frac{\sum_\ell b_\ell a(p^{\ell})}{p^s}\biggr)
      =\prod_{0\leq \ell \leq L}L(s,\sym^\ell f)^{b_\ell}H(s),
      \end{equation*}
      where
$H$ is a function that is  holomorphic and bounded by a  constant 
 in the region $\Re s\geq \frac12+\varepsilon$ for any $\ve >0$.
\end{lem}

\begin{proof}
This follows easily by comparing the $p$-th Euler factors.
\end{proof}

We recall that, in the half-plane of absolute convergence, we have
\begin{equation}
  \label{eq:3}
  L(s,f)=\prod_{p}
    \biggl(1-\frac{a(p)}{p^{s}}+\frac{1}{p^{2s}}\biggr)^{-1}
    =\prod_{p}
    \biggl(1-\frac{\alpha(p)}{p^{s}}\biggr)^{-1}\biggl(1-\frac{\beta(p)}{p^{s}}\biggr)^{-1}
\end{equation}
as well as
\begin{equation}
  \label{sym}
  L(s,\sym^2f)=\prod_{p}
    \biggl(1-\frac{a(p)^2-1}{p^{s}}+\frac{a(p)^2-1}{p^{2s}}-\frac{1}{p^{3s}}\biggr)^{-1}.
\end{equation}
\subsection{Averages of multiplicative functions}

We quote  Theorem 21.2 of \cite{Ramare*06} which follows an idea of Wirsing  \cite{Wirsing*61}.
\begin{lem}\label{density+}
  Let $f$ be a non-negative multiplicative function
and  $\kappa$ be a non-negative real parameter
such that
\begin{equation*}
  \left\{
    \begin{array}{l}
      \displaystyle
      \sum_{\substack{ p\geq2, \nu\geq1\\ p^{\nu}\leq Q}}
      f\bigl(p^{\nu}\bigr)\log\bigl(p^{\nu}\bigr)=
      \kappa Q+\Ocal(Q/\log(2Q))\qquad(Q\geq1),
      \\\displaystyle
      \sum_{p\geq2}
      \sum_{\substack{\nu,\ell\geq 1,\\ p^{\nu+\ell}\leq Q}}
      f\bigl(p^\ell\bigr)f\bigl(p^{\nu}\bigr)\log\bigl(p^{\nu}\bigr)
      \ll \sqrt{Q},
    \end{array}
  \right.
\end{equation*}
then we have
$$
\sum_{d\leq D}f(d)=  \kappa\, C\cdot D\left(\log D\right)^{\kappa-1} (1+o(1)),
$$
where
\begin{equation}
  \label{defC}
  C=\frac{1}{\Gamma(\kappa+1)}
      \prod_{p}\biggl\{
      \biggl(1-\frac1p\biggr)^{\kappa}
      \sum_{\nu\geq0}f\bigl(p^{\nu}\bigr)
      \biggr\}.
\end{equation}
\end{lem}
\begin{lem}
  \label{dens}
  Under the same hypotheses of Lemma~\ref{density+} we have, for any
  continuously differentiable function~$\eta$ with~$\int_0^1 (\log u)^\kappa \eta(u) du \neq 0$:
  \begin{equation*}
    \sum_{d\leq D}f(d)\eta(d/D)=\kappa C
    (1+o(1))\int_2^D (\log u)^{\kappa-1}\eta(u/D)du
  \end{equation*}
  provided the RHS goes to infinity as $D\longrightarrow \infty$.
\end{lem}
The condition on~$\eta$ is obviously satisfied if, as will be the case for us, $\eta$ is a non-negative with support inside the interval~$[0, 1]$.

\begin{proof}
  We find that
  \begin{align*}
    \sum_{d\leq D}f(d)\eta(d/D)
    &=\sum_{d\leq D}f(d)\eta(1)-\sum_{d\leq D}f(d)\int_{d/D}^1\eta'(t)dt
    \\&=\sum_{d\leq D}f(d)\eta(1)-\int_{0}^1\sum_{d\leq tD}f(d)\eta'(t)dt
    \\&=\kappa\eta(1) C\cdot D\left(\log D\right)^{\kappa-1}(1+o(1))
      \\&\qquad  -\int_{2}^D\kappa C u(\log
          u)^{\kappa-1}(1+o(1))\eta'(u/D)du/D
          +\Ocal(1)
  \end{align*}
  where the $\Ocal(1)$ is here to take care of the part of the
  integration between $1$ and $2$ since the $(\log u)^{\kappa-1}$ may
  give us some trouble.
  Hence
  \begin{equation*}
    \sum_{d\leq D}f(d)\eta(d/D)(\kappa C)^{-1}
    =
    \int_2^D (\kappa-1+\log u)(\log u)^{\kappa-2}(1+o(1))\eta(u/D)du
    +\Ocal(1)
  \end{equation*}
  from which we readily deduce that
  \begin{equation*}
    \sum_{d\leq D}f(d)\eta(d/D)(\kappa C)^{-1}
    =
    C (1+o(1)) D (\log D)^{\kappa - 1},
  \end{equation*}
  where~$C = \int_0^1 (\log(1/u))^{\kappa-1} \eta(u)du \neq 0$.
  On the other hand, we also have
  $$ \int_2^D (\log u)^{\kappa-1} \eta(u/D) du = C (1 + o(1)) D (\log D)^{\kappa - 1}, $$
  whence our claimed estimate.
\end{proof}

\section{A general average bound}

\begin{lem}
  \label{Usual}
  Let $L\in\N_{>0}$, and assume $(\myH_\ell)_{1\leq\ell\le L}$.
  Let $(b_\ell)_{0\leq \ell\leq L}$ be a collection of non-negative integers. Given a primitive form $f(z)=\sum_{n\geq 1} a(n)e(nz)$ as in the introduction, 
  let us define a multiplicative function $h_f$ by the equality
  \[
  \sum_n \frac{h_f(n)}{n^s}=\prod_{p}\biggl(1+\frac{\sum_\ell b_\ell a(p^{\ell})}{p^s}\biggr)
  \]
  Then $h_f$ is supported on square-free integers and there exists a polynomial $P_L$ of degree at most $b_0-1$ such that, for any $\varepsilon>0$, we have
  \begin{equation}\label{eq:equation-lemma-usual}
    \sum_{n\geq1 }
h_f(n)
    \eta(n/X)
    = XP_L(\log X)+\Ocal\Bigl(X^{\frac12+\varepsilon}\prod_{1\leq \ell \leq L}q(\sym^{\ell} (f))^{b_\ell\lambda_\ell+\ve}\Bigr)
  \end{equation}
  for any compactly supported twice continuously differentiable non-negative function $\eta$. 
\end{lem}

\begin{proof}
  Let us denote by $S$ the left-hand side of~\eqref{eq:equation-lemma-usual}.
  By taking Mellin transforms (\emph{e.g.} p.90 of~\cite{Iwaniec-Kowalski*04}), we get
\begin{equation*}
  S=\frac{1}{2i\pi}\int_{2-i\infty}^{2+i\infty}X^s \check{\eta}(s) ^s
  \sum_{n\geq 1}\frac{h_f(n)}{n^{s}}ds.
\end{equation*}
The fact that $\eta$ is twice continuously differentiable ensures us that its
Mellin transform verifies $\check{\eta}(s)\ll 1/(1+|s|^2)$ uniformly
in any closed vertical strip in the half plane $\Re s >0$.   
Lemma~\ref{genDir} gives us an expression for the Dirichlet series
$\sum_{n\geq1}h_f(n)/n^s$ from which we see that we can shift the line of
integration to $\Re s=\frac12+\varepsilon$ obtaining that the error term is at most
\begin{equation*}
  \Ocal\Bigl(X^{\frac12+\varepsilon}\prod_{1\leq \ell \leq L}q(\sym^{\ell} (f))^{b_\ell\lambda_\ell+\ve}\Bigr),
\end{equation*}
by our hypothesis $(\myH_\ell)_{1\leq\ell\le L}$ and the convexity principle. 
The residue at~$1$ gives the claimed main term, and the lemma follows readily.
\end{proof}

\section{A general Lemma around Vinogradov's trick}

\begin{lem}\label{Vin}
  Let $g$ be a real-valued multiplicative function supported on the
  squarefree integers. We assume further that $g(p)\geq F$ for every
  prime~$p$, and that for every prime $p\leq P$, we have $g(p)\geq \kappa>0$.
  Let $\eta$ be a non-negative, continuously
  differentiable function with support within $[0,1]$
  such that $\int_0^1\mkern-3mu\eta(v)dv=1$. 
  We have, for $M=P^\theta$ for some $\theta\in[0,1]$,
  \begin{equation*}
    \sum_{n\geq1}\mu^2(n)g(n)\eta\biggl(\frac{n}{PM}\biggr)
    \geq
    (1+o(1))\kappa\, C MP(\log
    MP)^{\kappa-1}
    \bigl(
    1-(\kappa-F)
    \myF(\theta;\kappa)\bigr)
  \end{equation*}
  where $C$ is given by~\eqref{defC} and $\myF$ is defined in~\eqref{defmyF}
\end{lem}
The factor $\mu^2(n)$ is only here to remind the reader that the
variable~$n$ is restricted to squarefree values. It can be omitted!
\begin{proof}
  We set
\begin{equation}
  \label{defS1}
  S=\sum_{n\geq1}g(n)\eta\biggl(\frac{n}{PM}\biggr).
\end{equation}
By our hypotheses, we find that
\begin{align*}
  S
  &=\sum_{\substack{n\leq PM,\\ P^+(n)\leq P}}
  g(n)\eta\biggl(\frac{n}{PM}\biggr)
  +
  \sum_{\substack{P<p\leq PM}}g(p)\sum_{n\leq PM/p}
  g(n)\eta\biggl(\frac{pn}{PM}\biggr)
  \\&\geq
  \sum_{\substack{n\leq PM,\\ P^+(n)\leq P}}
  g(n)\eta\biggl(\frac{n}{PM}\biggr)
  +
  F \sum_{\substack{P<p\leq PM}}\sum_{n\leq PM/p}
  \mu^2(n)\kappa^{\omega(n)}\eta\biggl(\frac{pn}{PM}\biggr)
  \\&\geq\sum_{\substack{n\leq PM}}
  \mu^2(n)\kappa^{\omega(n)}\eta\biggl(\frac{n}{PM}\biggr)
  +
  (F-\kappa)\sum_{\substack{P<p\leq PM}}\sum_{n\leq PM/p}
  \mu^2(n)\kappa^{\omega(n)}\eta\biggl(\frac{pn}{PM}\biggr).
\end{align*}
Here $P^+(n)$ denotes the greatest prime divisor of $n$.
We appeal to Lemma~\ref{dens} with $f(n)=\mu^2(n) \kappa^{\omega(n)}$ and get
\begin{multline*}
  S/(C \kappa)\geq
  (1+o(1))\int_2^{PM} (\log u)^{\kappa-1}\eta\Bigl(\frac{u}{PM}\Bigr)du
  \\
  +
  (F-\kappa+o(1))
  \sum_{\substack{N<p\leq PM}}\int_2^{PM/p} (\log u)^{\kappa-1}\eta\Bigl(\frac{up}{PM}\Bigr)du.
\end{multline*}
Note that the change of variable $vPM=u$ shows that
\begin{equation*}
  \int_2^{PM} (\log u)^{\kappa-1}\eta\Bigl(\frac{u}{PM}\Bigr)du
  =
  PM(\log PM)^{\kappa-1}\int_0^1\eta(v)dv(1+o(1)).
\end{equation*}
We use this estimate with $M$ replaced by $M/t$ and the prime number theorem to infer that
\begin{multline*}
  \sum_{\substack{N<p\leq PM}}\int_2^{PM/p} (\log
  u)^{\kappa-1}\eta\Bigl(\frac{up}{PM}\Bigr)du
  \\=
  PM(1+o(1))\int_0^1\eta(v)dv
  \int_N^{PM}
  \Bigr(\log\frac{PM}{t}\Bigl)^{\kappa-1}\frac {dt}{t\log t}
\end{multline*}
while this last integral equals, with the change of variable
$v=(PM)^h$ and $M=N^\theta$,
\begin{equation*}
  \int_1^{M}
  \frac {(\log v)^{\kappa-1}dv}{v(\log(PM)-\log v)}
  =
  (\log PM)^{\kappa-1}\int_0^{\theta/(1+\theta)}
  \frac {h^{\kappa-1}dh}{1-h}.
\end{equation*}
Recall that $\int_0^1\eta(v)dv=1$.
We thus find that
\begin{align*}
  \frac{(1+o(1))S}{C\kappa PM (\log
    PM)^{\kappa-1}}
  &\geq 
  1
  +(F-\kappa)\int_0^{\theta/(1+\theta)}
  \frac {h^{\kappa-1}dh}{1-h}\\
  &=1-(\kappa -F)\myF(\theta, \kappa).
\end{align*}
\end{proof}
\section{Proof of Theorems~\ref{geneV}, \ref{thm1}, \ref{thm2}, \ref{thmpos}}

Suppose $a(p) \not\in I$ for every $p \leq P$. Under the assumptions of Theorem~\ref{geneV}, let $\theta\in[0, 1]$ be such that
\begin{equation}\label{theta}
  \frac{1}{\kappa-F} > \myF(\theta;\kappa);
\end{equation}
for instance, we may take $\theta=\max(\myG(\kappa-F;\kappa)-\ve, 0)$.
Consider the sum 
\[
S=\sum_{n\ge1 }h_f (n)
    \eta(n/PM)
\]
where~$M\in[1, P]$. From  the upper and the lower bound of $S$ as given by Lemma~\ref{Usual} and~\ref{Vin} respectively and noting that $b_0=0$, we obtain,
\begin{equation*}
  (PM)^{\frac 1 2+\varepsilon}\prod_{1\leq \ell \leq L}q(\sym^{\ell} (f))^{b_\ell\lambda_{\ell}+\ve }
  \gg
  PM.
\end{equation*}
Therefore, with $M=P^ \theta$ for some $\theta\in[0,1]$ satisfying \eqref{theta}, we have
\begin{equation*}
  P \ll_k N^{\frac{2\sum_{\ell}\ell b_{\ell} \lambda_{\ell}}{1+\theta}+\ve}.
\end{equation*}
This leads to the estimate \eqref{vN} and the other estimate
\eqref{vk} is proved in a similar manner. 

Let us inspect what this gives to us under the convexity bound for
$\lambda_\ell=1/4$. Since the quantity $2\sum_{\ell\ge1}\ell 
b_\ell\lambda_\ell$ takes all the values that are half-positive integers, we
may inspect the first of them one by one. As we did above, we focus on the level~$N$.

\subsection*{First case $(1/2)\sum_{\ell\geq1} \ell b_\ell=1/2$}

This is only possible with the choice $b_1=1$, all other $b_\ell$'s
being~$0$. We have $\sum_{1\leq \ell\leq L}b_\ell U_\ell(x/2)=x$ which is
positive when $x=a(p)>0$. On assuming $a(p)\geq \delta$ when $p\leq P$,
we see that we may take $\kappa=\delta$ and $F=-2$ and get, for~$N\geq N_0(\ve)$,
\begin{equation}
  \label{eq:9.1}
  \frac{\log P}{\log N} \leq \frac{2\lambda_1}{1+\myG(2+\delta;\delta)}+\ve.
\end{equation}
Hence Theorem~\ref{thm1}.

\subsection*{Second case $(1/2)\sum_{\ell\geq1} \ell b_\ell=1$}

This is only possible with the choice $b_2=1$, all other $b_\ell$'s
being~$0$. We have $\sum_{1\leq \ell\leq L}b_\ell U_\ell(x/2)=x^2-1$ which is
positive when $x=a(p)\notin[-1,1]$. On assuming $|a(p)|\geq 1+\delta$ when $p\leq P$,
we see that we may take $\kappa=2\delta+\delta^2$ and $F=-1$ and get, for~$N\geq N_0(\ve)$,
\begin{equation}
  \label{eq:9.2}
  \frac{\log P}{\log N} \leq \frac{4\lambda_2}{1+\myG(1+2\delta+\delta^2;2\delta+\delta^2)}+\ve.
\end{equation}
Hence Theorem~\ref{thm2}.

\subsection*{Finding non-negative values}

Let~$I=[0, 2]$. A numerical computation found the coefficients~$(b_\ell)_{0\leq \ell \leq 5} = (0, 0, 3, 5, 4, 1)$, which satisfy~\eqref{eqI} with~$\kappa \geq 1/3$ and~$F = -10$. Then Theorem~\ref{thmpos} follows from the bounds~\eqref{vN} and~\eqref{vk}.

\section{Proof of Theorem \ref{thmnew}}

\medskip{}

Let~$\eta:\R_+ \to \R_+$ be smooth, compactly supported and such that~$\1_{[0, 1]} \geq \eta \geq \1_{[1/3, 2/3]}$, let~$\ve>0$, and consider
\[
T(X)=\sum_n {\mu}^2 (n) a(n) \eta(n/X) \qquad (X\geq 1).
\]
By Lemma~\ref{Usual}, we get
\begin{equation}\label{ub}
T(X) \ll X^{1/2}(k^2N)^{1/4+\ve}.
\end{equation}
Suppose that $a(p) \geq  0$ for all primes $p \leq X$. If the inequality
\begin{equation}
  \sum_{\substack{n:\\ a(n) \geq 1}} {\mu}^2 (n) a(n) \eta(n/X) \geq X^{1-\ve}\label{eq:case-distinction}
\end{equation}
holds then we easily have 
\begin{equation}\label{lb1}
T(X) \geq X^{1-\ve}.
\end{equation}
Otherwise, suppose that~\eqref{eq:case-distinction} does not hold. We write
 \begin{align*}
 T(X)&{}=\sum_{\substack{n:\\ 0 \leq a(n) <1}} {\mu}^2 (n) a(n) \eta(n/X) +
 \sum_{\substack{n:\\ a(n) \geq 1}} {\mu}^2 (n) a(n) \eta(n/X)\\
 &{}\geq \sum_{\substack{n:\\ 0 \leq a(n) <1}} {\mu}^2 (n) a(n)^2 \eta(n/X) +
 \sum_{\substack{n:\\ a(n) \geq 1}} {\mu}^2 (n) a(n) \eta(n/X)\\
 &{}=\sum_{n} {\mu}^2 (n) a(n)^2 \eta(n/X) +
 \sum_{\substack{n:\\ a(n) \geq 1}} {\mu}^2 (n) a(n)(1-a(n)) \eta(n/X).
 \end{align*}
 
Now the last sum is $\Ocal(X^{1-\ve/2})$ by Deligne's bound $|a(p)| \leq 2$ and the negation of~\eqref{eq:case-distinction}. The first sum can be handled by Rankin-Selberg method (Lemma~\ref{Usual}) and is $\gg L(1, \sym^2 f) X + \Ocal(X^{1/2}(k^2N^2)^{1/4+\ve})$.
Thus we have, using the lower bound $L(1, \sym^2 f) \gg 1/\log (kN)$ due to Hoffstein and Lockhart \cite{HL},
\begin{equation}\label{lb2}
T(X)\gg X/\log (kN) + \Ocal(X^{1/2}(k^2N^2)^{1/4+\ve}) +\Ocal (X^{1-\ve/2}).
\end{equation}

One of the equations~\eqref{lb1} and~\eqref{lb2} must hold and either, in conjunction with equation \eqref{ub}, imply the theorem.

\section{Proof of Theorem \ref{thm3}}

By equation~(3) of~\cite{Ramare*10}, the Deligne bound~$\abs{a(p)}\leq 2$ and Mertens' theorem (see \cite[Eq.~(2.15)]{Iwaniec-Kowalski*04}), we have
$$ \log L(1, f) = \Ocal_\ve(1) + \sum_{p\leq q^\ve} \frac{a(p)}p, $$
and therefore
$$ \sum_{p\leq q^\ve} \frac1p \Big(a(p) - \frac{\log L(1, f)}{\log\log q}\Big) = \Ocal_\ve(1). $$
However, if we had~$a(p)<\gamma^- - \delta$ for~$p\leq q^\ve$, then we would also have
$$ \sum_{p\leq q^\ve} \frac1p \Big(a(p) - \frac{\log L(1, f)}{\log\log q}\Big) \leq \Ocal_{\ve, \delta}(1) - \frac{\delta}2 \log\log q, $$
which is a contradiction for~$q$ large enough, and therefore there must be a prime~$p\leq q^\ve$ such that~$a(p)\geq \gamma_- - \delta$. An identical argument shows the existence of~$p\leq q^\ve$ such that~$a(p) \leq \gamma_+ + \delta$.

\section{Conditional bounds: proof of Theorem~\ref{th:ankeny}}

By the Stone-Weierstrass theorem, the fact that~$(U_\ell)$ forms a basis of~$\R[X]$, and the relation~\eqref{eq:5}, we may find~$L\geq 1$ and real
coefficients~$b_0, \dotsc, b_L$ depending on~$I$, with~$b_0>0$, such
that
\begin{equation}
  \sum_{p\leq x} \1(a(p)\in I) (1-\tfrac px) \log p \geq \sum_{\ell = 0}^L b_\ell \sum_{p\leq x} a_{\sym^\ell f} (p) (1-\tfrac px) \log p. \label{eq:condbound_mino}
\end{equation}
By Chebyshev's estimate, the contribution of the term~$\ell = 0$ is
$$ b_0 \sum_{p\leq x} (1-\tfrac px)\log p \gg_I x $$
with an absolute constant. To show that the right-hand side
of~\eqref{eq:condbound_mino} is positive for
some~$x = \Ocal_I((\log q)^2)$, it therefore suffices to show that for all
integer~$\ell \geq 1$ and all real~$x\geq 1$, we have
$$ \sum_{p\leq x} a_{\sym^\ell f}(p) (1-\tfrac px) \log p = \Ocal_\ell(x^{1/2} \log q). $$
This is an immediate consequence of the explicit formula
\cite[eq.~(5.33)]{Iwaniec-Kowalski*04} (with an additional smoothing,
as in~\cite[eq.~(13.28)]{Montgomery-Vaughan*06}) along with classical
zero density estimates~\cite[Theorem~5.8]{Iwaniec-Kowalski*04}.

\bibliographystyle{amsplain}

\begin{thebibliography}{10}

\bibitem{Ankeny*52}
N.~C. Ankeny, \emph{The least quadratic non residue}, Ann. of Math. (2)
  \textbf{55} (1952), 65--72.

\bibitem{Burgess*57}
D.~A. Burgess, \emph{The distribution of quadratic residues and non-residues},
  Mathematika \textbf{4} (1957), 106--112.

\bibitem{Clozel-Harris-Taylor*08}
L.~Clozel, M.~Harris, and R.~Taylor, \emph{Automorphy for some {$l$}-adic lifts
  of automorphic mod {$l$} {G}alois representations}, Inst. Hautes {\'E}tudes
  Sci. Publ. Math. (2008), no.~108, 1--181.

\bibitem{ClozelThorne2015}
L.~Clozel and J.~A. Thorne, \emph{Level raising and symmetric power
  functoriality, {II}}, Ann. of Math. (2) \textbf{181} (2015), no.~1, 303--359.

\bibitem{ClozelThorne2017}
\bysame, \emph{Level-raising and symmetric power functoriality, {III}}, Duke
  Math. J. \textbf{166} (2017), no.~2, 325--402.

\bibitem{Deligne1974}
P.~Deligne, \emph{La conjecture de {W}eil. {I}}, Inst. Hautes {\'E}tudes Sci.
  Publ. Math. (1974), no.~43, 273--307.

\bibitem{GelbartJacquet1978}
S.~Gelbart and H.~Jacquet, \emph{A relation between automorphic
  representations of {${\rm GL}(2)$} and {${\rm GL}(3)$}}, Ann. Sci. \'{E}cole
  Norm. Sup. (4) \textbf{11} (1978), no.~4, 471--542.

\bibitem{GodementJacquet}
R.~Godement and H.~Jacquet, \emph{Zeta functions of simple algebras}, Lecture
  Notes in Mathematics, Vol. 260, Springer-Verlag, Berlin-New York, 1972.
  



\bibitem{HarrisEtAl2010}
M.~Harris, N.~Shepherd-Barron, and R.~Taylor, \emph{A family of {C}alabi-{Y}au
  varieties and potential automorphy}, Ann. of Math. (2) \textbf{171} (2010),
  no.~2, 779--813.
  
\bibitem{HL}
J.~ Hoffstein and P.~Lockhart, \emph{Coefficients of Maass forms and the Siegel zero. With an appendix by Dorian Goldfeld, Hoffstein and Daniel Lieman.}
 Ann. of Math. (2) 140 (1994), no. 1, 161--181.
  
  \bibitem{HR}
J.~Hoffstein, D.~Ramakrishnan, \emph{Siegel zeros and cusp forms}, Int. Math. Res. Not., Volume 1995, Issue 6, 1995, Pages 279–308,
1995.

\bibitem{Iwaniec1997}
H.~Iwaniec, \emph{Topics in classical automorphic forms}, Graduate {Studies} in
  {Mathematics}, vol.~17, American Mathematical Society, Providence, RI, 1997.

\bibitem{Iwaniec-Kohnen-Sengupta*07}
H.~Iwaniec, W.~Kohnen, and J.~Sengupta, \emph{The first negative {H}ecke
  eigenvalue}, Int. J. Number Theory \textbf{3} (2007), no.~3, 355--363.


\bibitem{Iwaniec-Kowalski*04}
H.~Iwaniec and E.~Kowalski, \emph{Analytic number theory}, American
  Mathematical Society Colloquium Publications, American Mathematical Society,
  Providence, RI, 2004, xii+615 pp.

\bibitem{JM}
M. Jutila and Y. Motohashi, 
  \emph{Uniform bounds for Hecke $L$-functions}, 
  Acta Math., Vol. 195, No. 1, 2005, pp. 61--115

\bibitem{KimShahidi2002a}
H.~H. Kim and F.~Shahidi, \emph{Cuspidality of symmetric powers with
  applications}, Duke Math. J. \textbf{112} (2002), no.~1, 177--197.

\bibitem{KimShahidi2002b}
\bysame, \emph{Functorial products for {${\rm GL}_2\times{\rm GL}_3$} and the
  symmetric cube for {${\rm GL}_2$}}, Ann. of Math. (2) \textbf{155} (2002),
  no.~3, 837--893, With an appendix by Colin J. Bushnell and Guy Henniart.

\bibitem{Kim2003}
H.~H. Kim, \emph{Functoriality for the exterior square of {${\rm GL}_4$} and
  the symmetric fourth of {${\rm GL}_2$}}, J. Amer. Math. Soc. \textbf{16}
  (2003), no.~1, 139--183, With appendix 1 by Dinakar Ramakrishnan and appendix
  2 by Kim and Peter Sarnak.

\bibitem{Kowalski-Lau-Soundararajan-Wu*10}
E.~Kowalski, Y.-K. Lau, K.~Soundararajan, and J.~Wu, \emph{On modular signs},
  Math. Proc. Cambridge Philos. Soc. \textbf{149} (2010), no.~3, 389--411.
  \MR{2726725}

\bibitem{Langlands} R.P. Langlands, \emph{Problems in the Theory of Automorphic Forms}, Lecture
Notes in Mathematics, Vol. 170, Springer-Verlag, New York, 1970, pp.
18-86.

\bibitem{Lau-Liu-Wu*12}
Y. K. Lau, J. Y. Liu, and J.~Wu, \emph{The first negative coefficients of
  symmetric square {$L$}-functions}, Ramanujan J. \textbf{27} (2012), no.~3,
  419--441.

\bibitem{LemkeOliver-Thorner*07}
R.~J. Lemke~Oliver and J. Thorner, \emph{Effective log-free zero density
  estimates for automorphic {L}-functions and the {Sato}-{Tate} conjecture},
  Internat. Math. Res. Notices (2018), 30pp.
  
\bibitem{LW}
Y.-K. Lau and J. Wu, \emph{A density theorem on automorphic L-functions and some applications}, Trans. Amer. Math. Soc. 358 (2006), 441--472.

\bibitem{LW2}
Y.-K. Lau and J. Wu, \emph{Extreme values of symmetric power L-functions at 1}, Acta Arith.
126 (2007), No. 1, 57--76.

\bibitem{Matomaki*12}
K. Matom\"aki, \emph{On signs of {F}ourier coefficients of cusp forms},
  Math. Proc. Cambridge Philos. Soc. \textbf{152} (2012), no.~2, 207--222.
  \MR{2887873}

\bibitem{Montgomery-Vaughan*06}
H.L. Montgomery and R.C. Vaughan, \emph{Multiplicative {N}umber {T}heory: {I}.
  {C}lassical {T}heory}, Cambridge Studies in Advanced Mathematics, vol.~97,
  Cambridge University Press, 2006.
  
\bibitem{Munshi}
R. Munshi, \emph{The subconvexity problem for $L$-functions},
  Proc. Int. Cong. of Math. -- 2018
Rio de Janeiro, Vol. 1 (361--374), available at
https://eta.impa.br/dl/164.pdf


\bibitem{NewtonThorne2020}
J. Newton, J. A. Thorne : \emph{Symmetric power functoriality for holomorphic modular forms},
  Publ. Math. Inst. Hautes {\'E}tudes Sci. \textbf{134} (2021), 1--116, https://doi.org/10.1007/s10240-021-00127-3.

\bibitem{Ramare*06}
O.~Ramar{\'e}, \emph{Arithmetical aspects of the large sieve inequality},
  Harish-Chandra Research Institute Lecture Notes, vol.~1, Hindustan Book
  Agency, New Delhi, 2009, With the collaboration of D. S. Ramana.
  \MR{MR2493924}

\bibitem{Ramare*10}
O.~Ramar{\'e}, \emph{Comparing~$L(s, \chi)$ with its tructated Euler product and generalization},
  Funct. Approximatio, Comment. Math. \textbf{42} (2010), No. 2, 145--151.

\bibitem{Taylor*08}
R.~Taylor, \emph{Automorphy for some {$l$}-adic lifts of automorphic mod {$l$}
  {G}alois representations. {II}}, Publ. Math. Inst. Hautes {\'E}tudes Sci.
  (2008), no.~108, 183--239.

\bibitem{Vinogradov*27b}
I.~M. Vinogradov, \emph{On a general theorem concerning the distribution of the
  residues and non-residues of powers}, Trans. Amer. Math. Soc. \textbf{29}
  (1927), no.~1, 209--217.

\bibitem{Wirsing*61}
{E.} Wirsing, \emph{Das asymptotische {V}erhalten von {S}ummen {\"u}ber
  multiplikative {F}unktionen}, Math. Ann. \textbf{143} (1961), 75--102.

\end{thebibliography}

\providecommand{\bysame}{\leavevmode\hbox to3em{\hrulefill}\thinspace}
\providecommand{\MR}{\relax\ifhmode\unskip\space\fi MR }
\providecommand{\MRhref}[2]{%
  \href{http://www.ams.org/mathscinet-getitem?mr=#1}{#2}
}
\providecommand{\href}[2]{#2}

\end{document}